\documentclass[12pt, english]{article}

\usepackage[english]{babel}
\usepackage[applemac]{inputenc}
\usepackage[T1]{fontenc} 

\usepackage{pgf}
\usepackage{tikz}
\usetikzlibrary{arrows}
\usetikzlibrary{decorations.pathmorphing}
\usetikzlibrary{backgrounds}
\usetikzlibrary{fit}

\usepackage{amsmath}
\usepackage{amssymb}
\usepackage{amsthm}
\usepackage{url}

\usepackage{multicol}

\usepackage{chessboard}
\usepackage{skak}
\usepackage{chessfss}
\usepackage{adjustbox}

\usepackage[LSBC1,T1]{fontenc}
\setboardfontencoding{LSBC1}\setboardfontsize{20pt}

\theoremstyle{plain}

\theoremstyle{definition}

\begin{document}

\title{The bishop and knight checkmate on a large chessboard}

\author{Johan Wästlund}

\date{} 

\maketitle

\begin{abstract}
In 1983 the chess periodical EG published a summary of a letter from Julius Telesin outlining how a king, a bishop and a knight can checkmate a lonely king on an arbitrarily large chessboard. The Telesin checkmating procedure doesn't seem to be widely known, and the published account left out a number of details.
We describe the method precisely and show that it works against every defense. We also discuss the open question of the asymptotics of the largest distance to mate on a large square board. 
\end{abstract}

\section{Introduction}
The purpose of this article is to demonstrate rigorously a fact that doesn't seem to be widely known and that hasn't been explained in much detail in the literature: A king, a bishop and a knight can force checkmate against a lonely king on an arbitrarily large chessboard! 

For this to be possible, the board must have a corner of the color that the bishop can reach. But provided this is the case, and excluding situations where there is an immediate problem of stalemate or a piece being lost, checkmate can always be achieved. Naturally the 50-move rule must be waived throughout this discussion.

The checkmating procedure was described by Julius Telesin in a letter summarized in the journal EG (a journal devoted to studies and theoretical aspects of the endgame) in 1983 \cite{Telesin}. There is no doubt about the correctness, but in some aspects the exposition is not very clear. 

More recently the question has been discussed at Chess Stack Exchange \cite{CSE} without a conclusive answer, and apparently without anyone being aware of Telesin's solution.    
In the thread there is an answer from 2015 by Dan Stronger, who gives data from a tablebase program showing that mate can be achieved in at most 93 moves on a 16 by 16 board. He says that (presumably by playing around with the program) he has become convinced that checkmate is possible on arbitrarily large boards, but only briefly describes the strategy.
Stronger's longest winning line is cited in the Wikipedia article \cite{Wikipedia} without mention of Telesin's earlier result. 

V\'aclav Kot\v{e}\v{s}ovec pushed the tablebase calculations to board size $24 \times 24$ in 2017, and gives the longest distance to mate in \cite[p.~349]{Kotesovec}, see also \cite{oeis}, linking to the Telesin article. Another interesting paper is \cite{GMSB}, which discusses the extraction of useful strategies from tablebases for the standard board. The KBN vs K checkmate on large boards has also been briefly discussed on MatPlus.Net \cite{MatPlus}. 

Telesin's solution (as summarized in \cite{Telesin}, presumably by an editor) is somewhat difficult to verify.
In some cases it is explained that, for instance, the knight ``approaches'', but before understanding the overall plan in enough detail, it's not entirely clear where exactly the knight wants to go. Neither is it obvious that all defenses have been considered. 
I therefore want to share my thoughts on the KBN vs K checkmate, making some points more precise. 

It may be argued that to the ordinary chess player, mating procedures on boards larger than 8 by 8 are of limited interest. But the reason that many players find the KBN vs K checkmate rather tricky might be precisely that one doesn't see how it would work in general. It appears ad-hoc, and therefore becomes hard to understand and remember. Some interesting similar remarks are made in \cite[p.168]{Borovik}. 

\section{Background}
We first summarize the standard chess knowledge about the KBN vs K endgame. We assume throughout that White has the two pieces and Black the lonely king. To force checkmate, the black king has to be driven into a corner that the white bishop can reach (otherwise there will be a stalemate problem). 

\begin{figure}[h]
\begin{center}
\begin{tikzpicture} [scale=0.75]


\foreach \x in { 4, 6}
  \foreach \y in {0, 2, 4}  
    \filldraw [ color=white!70!gray] (\x, \y) rectangle (\x+1, \y+1);
\foreach \x in {3, 5, 7}
  \foreach \y in {1, 3}  
    \filldraw [ color=white!70!gray] (\x, \y) rectangle (\x+1, \y+1);
\foreach \y in {0, 2, 4}
  \filldraw [color=white!70!gray] (2.5, \y) rectangle (3, \y+1);
\foreach \x in {3, 5, 7}
  \filldraw [color=white!70!gray] (\x, 5) rectangle (\x+1, 5.5);

\node  at (6.5, 1.45) {\BlackKingOnWhite};
\node  at (4.5, 3.45) {\WhiteKnightOnWhite};
\node  at (6.5, 3.45) {\WhiteBishopOnWhite};

\draw [ultra thick, color=red] (8, 2)--(6,2)--(6,1)--(4,1)--(4,0); 
\draw[thick] (2.5,0)--(8,0)--(8,5.5);
 
\node [anchor=base] at (3.5, -0.7) {d};
\node [anchor=base] at (4.5, -0.7) {e};
\node [anchor=base] at (5.5, -0.7) {f};
\node [anchor=base] at (6.5, -0.7) {g};
\node [anchor=base] at (7.5, -0.7) {h};

\node at (8.5, 0.5) {$1$};
\node at (8.5, 1.5) {$2$};
\node at (8.5, 2.5) {$3$};
\node at (8.5, 3.5) {$4$};
\node at (8.5, 4.5) {$5$};

\begin{scope}[xshift=8.5cm]


\foreach \x in { 4, 6}
  \foreach \y in {0, 2, 4}  
    \filldraw [ color=white!70!gray] (\x, \y) rectangle (\x+1, \y+1);
\foreach \x in {3, 5, 7}
  \foreach \y in {1, 3}  
    \filldraw [ color=white!70!gray] (\x, \y) rectangle (\x+1, \y+1);
\foreach \y in {0, 2, 4}
  \filldraw [color=white!70!gray] (2.5, \y) rectangle (3, \y+1);
\foreach \x in {3, 5, 7}
  \filldraw [color=white!70!gray] (\x, 5) rectangle (\x+1, 5.5);

\node  at (7.5, 1.45) {\BlackKingOnWhite};
\node  at (5.5, 2.45) {\WhiteKnightOnWhite};
\node  at (6.5, 1.45) {\WhiteBishopOnWhite};
\node at (5.5, 1.45) {\WhiteKingOnWhite};

\draw[thick] (2.5,0)--(8,0)--(8,5.5);
 
\node [anchor=base] at (3.5, -0.7) {d};
\node [anchor=base] at (4.5, -0.7) {e};
\node [anchor=base] at (5.5, -0.7) {f};
\node [anchor=base] at (6.5, -0.7) {g};
\node [anchor=base] at (7.5, -0.7) {h};

\node at (8.5, 0.5) {$1$};
\node at (8.5, 1.5) {$2$};
\node at (8.5, 2.5) {$3$};
\node at (8.5, 3.5) {$4$};
\node at (8.5, 4.5) {$5$};

\end{scope}

\end{tikzpicture}
\caption{Left: Once the black king is constrained to the immediate vicinity of a corner that the bishop can reach, checkmating is relatively straightforward. Right: The white king approaches, and Black is finally checkmated either in the corner or on a square next to the corner.}
\label{F:finish}
\end{center}
\end{figure}
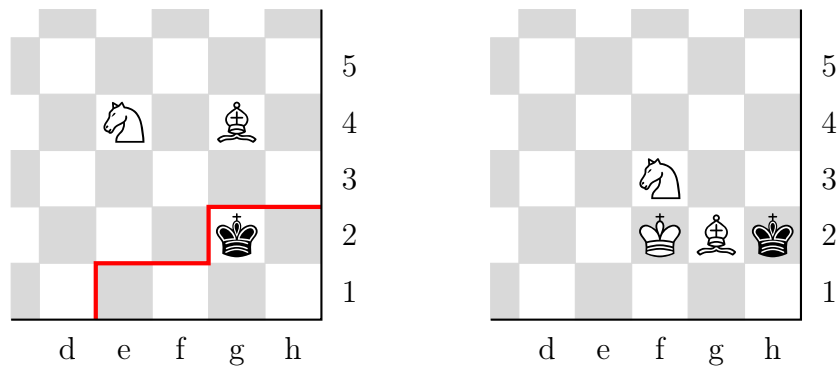

Once we obtain a situation as in Figure~\ref{F:finish}, the finish is relatively easy. Here the white pieces constrain the black king to a set of six squares near the corner. White can now maneuver their king to e2, then (in some order) the bishop via h3 to f1 and the king to f2. Finally the knight reroutes via g5 for instance, and checkmate is given by Bg2 and Nf3 in the correct order, giving check in the penultimate move.

On the standard 8 by 8 board, once White has centralized their own king, the black king is already fairly close to the edge of the board. Knowing that only two of the corners allow checkmate, the defending player often heads for one of the ``safe'' corners. The tricky part is then to force the black king over to another corner. This is what seems ad-hoc and leaves the impression that the method just barely works. 

For instance, the textbook instructions generated by \cite{GMSB} explicitly state that step 1 of the procedure is to drive the black king to the edge of the board, and step 2 to force the king to the appropriate corner. As we shall see, step 2 is achievable, but on a large board it seems we cannot drive the black king from a safe corner to another corner along the edge, and therefore it's not clear whether step 1 brings us any closer to achieving step 2. 

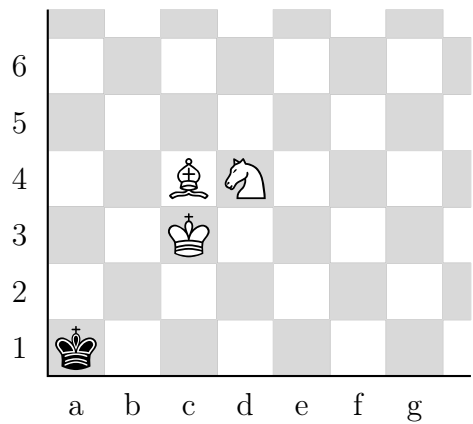
\begin{figure}[h]

\begin{center}

\begin{tikzpicture} [scale=0.75]

\draw[step=1cm,gray,ultra thin] (0,0) grid (7.5, 6.5); 

\foreach \x in {0, 2, 4, 6}
  \foreach \y in {0, 2, 4}  
    \filldraw [ color=white!70!gray] (\x, \y) rectangle (\x+1, \y+1);
\foreach \x in {1, 3, 5}
  \foreach \y in {1, 3, 5}  
    \filldraw [ color=white!70!gray] (\x, \y) rectangle (\x+1, \y+1);
\foreach \y in {1, 3, 5}
  \filldraw [color=white!70!gray] (7, \y) rectangle (7.5, \y+1);
\foreach \x in {0, 2, 4, 6}
  \filldraw [color=white!70!gray] (\x, 6) rectangle (\x+1, 6.5);

\node at (2.5,2.45) {\WhiteKingOnWhite};
\node at (0.5,0.45) {\BlackKingOnWhite};
\node  at (3.5, 3.45) {\WhiteKnightOnWhite};
\node  at (2.5, 3.45) {\WhiteBishopOnWhite};

\node [anchor=base] at (0.5, -0.7) {a};
\node [anchor=base] at (1.5, -0.7) {b};
\node [anchor=base] at (2.5, -0.7) {c};
\node [anchor=base] at (3.5, -0.7) {d};
\node [anchor=base] at (4.5, -0.7) {e};
\node [anchor=base] at (5.5, -0.7) {f};
\node [anchor=base] at (6.5, -0.7) {g};

\node at (-0.5, 0.5) {$1$};
\node at (-0.5, 1.5) {$2$};
\node at (-0.5, 2.5) {$3$};
\node at (-0.5, 3.5) {$4$};
\node at (-0.5, 4.5) {$5$};
\node at (-0.5, 5.5) {$6$};

\draw[thick] (0,6.5)--(0,0)--(7.5,0);

\end{tikzpicture}

\caption{White has to force the black king out of a ``safe'' corner.}
\label{F:digOut}
\end{center}
\end{figure}

Suppose as in Figure~\ref{F:digOut} that we have a light-squared bishop and want to drive the king from a1 to h1, or in case of a larger board, to a more distant light square corner on the first rank. To dig the black king out of the corner, the knight has to control a1. And if we get the black king to c1, some white piece must control b1 to stop Black from going back. Whether we play Na3, Ba2, or Ka2, some white piece will be misplaced. 

On an 8 by 8 board, we can still shove the king more or less directly from one corner to another, since the opposite edge stops the black king from running further: After 1.~Nc2+ Kb1 2.~Bb3 (waiting) Kc1 3.~Ba2 Kd1 4.~Nd4 Ke1 5.~Kd3 Kf2, we have 6.~Ne2! (this move has famously been missed even by strong grandmasters). In combination with Be6 (stopping by on b3 if Black tries to go back), this prevents the black king from escaping.

But on a much larger board, this would not work: 6.~- Kg2 7.~Be6 Kh2 followed by 8.~- Ki3, and Black runs towards the north-east corner. 

As was pointed out by a reviewer, there is another method of pushing the black king to the right corner known as the Delétang triangles method \cite{Deletang}. This method consists in restricting the black king diagonal by diagonal to ever smaller triangles, which is also how the Telesin method works. But it requires the knight to hold part of the barrier and the white king to move into the triangle, and this doesn't seem to generalize to a large board. 

\section{Telesin's method}
It seems that on a large board, White cannot force the black king from one corner to another along the edge of the board. But Telesin shows, modifying the triangles method, that we can constrain the black king not along the edge, but diagonal by diagonal. 

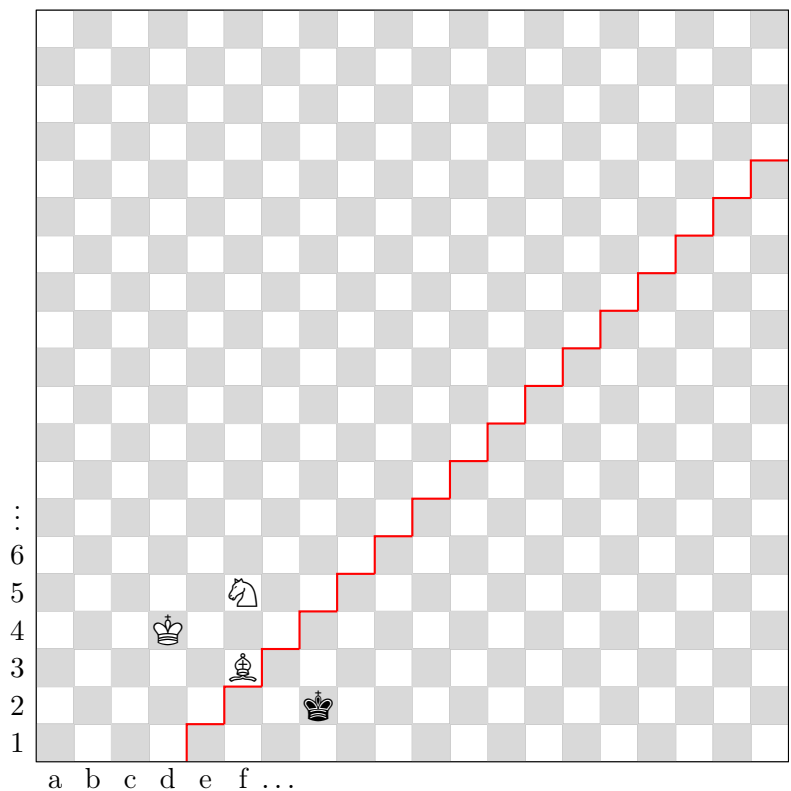
\begin{figure}[h]
\begin{center}

\begin{tikzpicture} [scale=0.5]

\draw[step=1cm,gray,ultra thin] (0,0) grid (20, 20); 

\foreach \x in {0, 2, 4, 6, 8, 10, 12, 14, 16, 18}
  \foreach \y in {0, 2, 4, 6, 8, 10, 12, 14, 16, 18}  
    \filldraw [ color=white!70!gray] (\x, \y) rectangle (\x+1, \y+1);
\foreach \x in {1, 3, 5, 7, 9, 11, 13, 15, 17, 19}
  \foreach \y in {1, 3, 5, 7, 9, 11, 13, 15, 17, 19}  
    \filldraw [ color=white!70!gray] (\x, \y) rectangle (\x+1, \y+1);

\node at (3.5,3.45) {\setboardfontsize{14pt}\WhiteKingOnWhite};
\node at (7.5,1.45) {\setboardfontsize{14pt}\BlackKingOnWhite};
\node  at (5.5, 4.45) {\setboardfontsize{14pt}\WhiteKnightOnWhite};
\node  at (5.5, 2.45) {\setboardfontsize{14pt}\WhiteBishopOnWhite};

\node [anchor=base] at (0.5, -0.7) {\small a};
\node [anchor=base] at (1.5, -0.7) {\small b};
\node [anchor=base] at (2.5, -0.7) {\small c};
\node [anchor=base] at (3.5, -0.7) {\small d};
\node [anchor=base] at (4.5, -0.7) {\small e};
\node [anchor=base] at (5.5, -0.7) {\small f};
\node [anchor=base] at (6.5, -0.7) {\small $\dots$};

\node at (-0.5, 0.5) {\small $1$};
\node at (-0.5, 1.5) {\small $2$};
\node at (-0.5, 2.5) {\small $3$};
\node at (-0.5, 3.5) {\small $4$};
\node at (-0.5, 4.5) {\small $5$};
\node at (-0.5, 5.5) {\small $6$};
\node at (-0.5, 6.7) {\small $\vdots$};

\foreach \x in {0,1,2,3,4,5,6,7,8,9,10,11,12,13,14,15}
\draw[thick, red] (\x+4,\x)--(\x+4,\x+1)--(\x+5,\x+1);

\draw (0,0)--(0,20)--(20,20)--(20,0)--cycle;
\end{tikzpicture}

\caption{Telesin's method on a large chess board.}
\label{F:large}

\end{center}
\end{figure}

We assume for simplicity that the board is square with an even number of squares along each side, although everything should work in the same way on any rectangular board with at least one even side. 
Suppose we start from the position of Figure~\ref{F:digOut} but on a large board. Still numbering from a dark a1 corner, Telesin \cite{Telesin} instead suggests (after 1.~Nc2+ Kb1 2.~Bb3 Kc1 3.~Ba2 Kd1 4.~Nd4 Ke1) the line 5.~Bd5 Kf2 6.~Nf5 Kg1 7.~Kd4 Kh2 8.~Bf3, leading to the position of Figure~\ref{F:large}.

It seems we have let the black king out, but as Telesin shows, White can restrict Black to the south-east side of a diagonal barrier as indicated. Black can run almost all the way up to the north-east corner, but must at some point turn back and allow White to decrease the triangular zone. 

\subsection*{The simple enclosure}

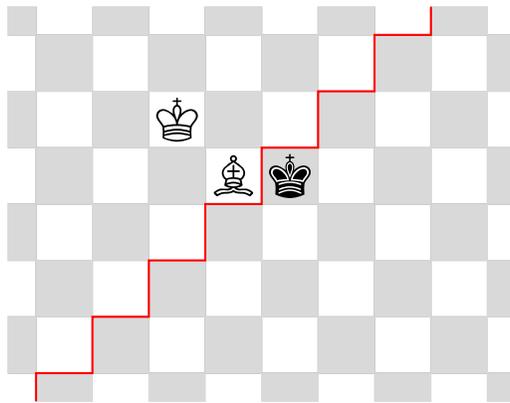
\begin{figure}[h]

\begin{center}
\begin{tikzpicture} [scale=0.75]

\draw[step=1cm, gray, ultra thin] (-3.5, 0) grid (5.5,6.5); 

\foreach \x in {-2, 0, 2, 4}
  \foreach \y in {0, 2, 4}  
    \filldraw [ color=white!70!gray] (\x, \y) rectangle (\x+1, \y+1);
\foreach \x in {-3, -1, 1, 3}
  \foreach \y in {1, 3,  5}  
    \filldraw [ color=white!70!gray] (\x, \y) rectangle (\x+1, \y+1);
\foreach \y in {0, 2, 4}
  \filldraw [color=white!70!gray] (-3.5, \y) rectangle (-3, \y+1);
\foreach \y in {1, 3, 5}
  \filldraw [color=white!70!gray] (5, \y) rectangle (5.5, \y+1);
\foreach \x in {-2, 0, 2, 4}
  \filldraw [color=white!70!gray] (\x, 6) rectangle (\x+1, 6.5);
\foreach \x in {-3, -1, 1, 3}
  \filldraw [color=white!70!gray] (\x, -0.5) rectangle (\x+1, 0);

\filldraw  [color=white!70!gray] (-3.5, 6) rectangle (-3, 6.5);
\filldraw  [color=white!70!gray] (5, -0.5) rectangle (5.5, 0);

\node  at (1.5, 3.45) {\BlackKingOnWhite};
\node  at (-0.5, 4.45) {\WhiteKingOnWhite};
\node  at (0.5, 3.45) {\WhiteBishopOnWhite};

\draw [thick, color=red] (4, 6.5)--(4,6)--(3,6)--(3,5)--(2,5)--(2,4)--(1,4)--(1,3)--(0,3)--(0,2)
--(-1,2)--(-1,1)--(-2,1)--(-2,0)--(-3,0)--(-3,-0.5);
 
\end{tikzpicture}

\caption{The simple enclosure}
\label{F:simple}
\end{center}

\end{figure}

The basic idea of the checkmating procedure is that the white king and bishop can keep the black king constantly on one side of a certain diagonal. We begin by demonstrating how this works, formalizing Telesin's concept of ``zone restriction'' \cite{Telesin}, and later turn to how to make progress and actually drive the black king towards a corner. Suppose that the kings and the white bishop are in the position of Figure~\ref{F:simple}.

In this situation, the white king and bishop can keep the black king on the south-east side of the indicated barrier. This can be done without moving the bishop, and without the help of the knight. White just has to follow Black along the barrier, making sure to always put the king on a file at least as far to the west as Black's king, and on a rank at least as far north, while also being in time to defend the bishop if Black attacks it. In order for this to work, the knight should be on the north-west side of the barrier so as not to allow Black to gain time by attacking it, and it should also be a couple of squares away from the barrier in order not to stand in the way of the king or bishop. 

\subsection*{The dented enclosure}

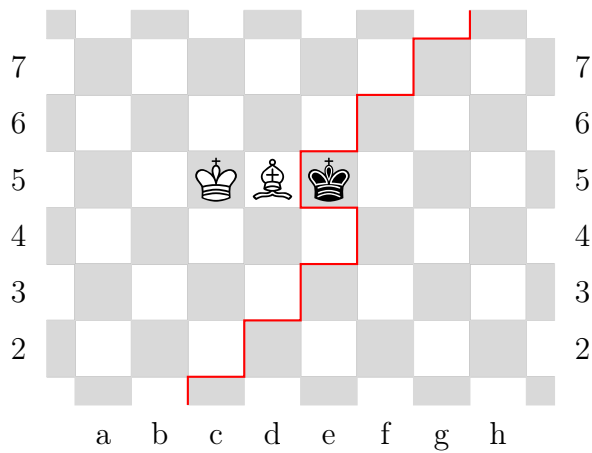
\begin{figure}[h]
\begin{center}
\begin{tikzpicture} [scale=0.75]

\draw[step=1cm, gray, ultra thin] (-3.5, 0) grid (5.5,6.5); 

\foreach \x in {-2, 0, 2, 4}
  \foreach \y in {0, 2, 4}  
    \filldraw [ color=white!70!gray] (\x, \y) rectangle (\x+1, \y+1);
\foreach \x in {-3, -1, 1, 3}
  \foreach \y in {1, 3,  5}  
    \filldraw [ color=white!70!gray] (\x, \y) rectangle (\x+1, \y+1);
\foreach \y in {0, 2, 4}
  \filldraw [color=white!70!gray] (-3.5, \y) rectangle (-3, \y+1);
\foreach \y in {1, 3, 5}
  \filldraw [color=white!70!gray] (5, \y) rectangle (5.5, \y+1);
\foreach \x in {-2, 0, 2, 4}
  \filldraw [color=white!70!gray] (\x, 6) rectangle (\x+1, 6.5);
\foreach \x in {-3, -1, 1, 3}
  \filldraw [color=white!70!gray] (\x, -0.5) rectangle (\x+1, 0);

\filldraw  [color=white!70!gray] (-3.5, 6) rectangle (-3, 6.5);
\filldraw  [color=white!70!gray] (5, -0.5) rectangle (5.5, 0);

\node  at (1.5, 3.45) {\BlackKingOnWhite};
\node  at (-.5, 3.45) {\WhiteKingOnWhite};
\node  at (0.5, 3.45) {\WhiteBishopOnWhite};

\draw [thick, color=red] (4, 6.5)--(4,6)--(3,6)--(3,5)--(2,5)--(2,4)--(1,4)--(1,3)--(2,3)--(2,2)
--(1,2)--(1,1)--(0,1)--(0,0)--(-1,0)--(-1,-0.5);
 
\node [anchor=base] at (-2.5, -1.2) {a};
\node [anchor=base] at (-1.5, -1.2) {b};
\node [anchor=base] at (-0.5, -1.2) {c};
\node [anchor=base] at (0.5, -1.2) {d};
\node [anchor=base] at (1.5, -1.2) {e};
\node [anchor=base] at (2.5, -1.2) {f};
\node [anchor=base] at (3.5, -1.2) {g};
\node [anchor=base] at (4.5, -1.2) {h};

\node at (-4, 0.5) {$2$};
\node at (-4, 1.5) {$3$};
\node at (-4, 2.5) {$4$};
\node at (-4, 3.5) {$5$};
\node at (-4, 4.5) {$6$};
\node at (-4, 5.5) {$7$};

\node at (6, 0.5) {$2$};
\node at (6, 1.5) {$3$};
\node at (6, 2.5) {$4$};
\node at (6, 3.5) {$5$};
\node at (6, 4.5) {$6$};
\node at (6, 5.5) {$7$};
 
\end{tikzpicture}

\caption{The dented enclosure.}
\label{F:dented}

\end{center}
\end{figure}

In Figure~\ref{F:simple}, as soon as Black makes a move, White can put their king on a dark square right next to the bishop and on the opposite side to the black king with respect to the bishop's north-west to south-east diagonal. White can then maintain what we will call a \emph{dented enclosure} (Figure~\ref{F:dented}), even further restricting Black's movement. Here we number the ranks and files without any assumptions on how far away the corners are. 

By applying a certain strategy, White can make sure that Black never crosses this ``dented'' barrier. If Black tries to escape to the north-east, White follows along with the king just as when defending the simple enclosure. If on the other hand Black goes via f4 to e3, White plays Kc4 creating what Telesin calls a ``valve'' (more about that in a moment), and unless Black goes back, White can install the bishop on e4, reducing the problem of containing the black king to a smaller simple enclosure.

An important point is that we can't allow Black to obtain a repetition of moves. If we are constantly kept busy, we can't make progress towards checkmate. 
If in Figure~\ref{F:dented} Black plays 1.~- Kf4, we can't answer 2.~Kd4, because on 2.~- Kf5 we would have to go back with 3.~Kc5 to be in time if Black tries Kf5-f6-e7. Since Black could then repeat with 3.~- Kf4, it's clear that we can't make progress by meeting Kf4 with Kd4. Similarly we can't meet Kf4 with a bishop move, since the bishop too would have to return on the following move. Therefore on 1.~Kf4, the king and bishop must wait, and the knight make a move towards its destination (which we discuss in a moment).

After 1.~- Kf4 2.~Nxx Ke3, the only response that defends the barrier is 3.~Kc4. And if Black then continues 3.~- Kd2, we have to play 4.~Be4 or the black king will get out. We see how several of the moves, including the ``valve'', are forced by the outline of the barrier. These moves are given in \cite{Telesin}, but become harder to understand without an exact description of the zone to which the black king is restricted.

If we get to play 4.~Be4, we have established a simple enclosure which is smaller than the dented enclosure we just had. And as soon as we have decreased the enclosure, we have made progress. Notice also that if Black runs to the north-east and White follows, then as soon as Black makes a move that doesn't go north, White can establish a smaller dented enclosure by moving the bishop along the barrier, to e6 or higher.

\subsection*{The knight's action if Black waits}

We have seen that a simple enclosure always leads to a smaller dented enclosure. Moreover, from a dented enclosure, as soon as the black king tries to run anywhere, the white king and bishop can decrease the enclosure one way or another. 
 It's only if Black lingers near the white king and bishop that the knight will have to take action. Therefore it doesn't matter how far away the knight is.  

\begin{figure}   
   
\begin{center}
\begin{tikzpicture} [scale=0.75]

\draw[step=1cm, gray, ultra thin] (-3.5, 0) grid (5.5,6.5); 

\foreach \x in {-2, 0, 2, 4}
  \foreach \y in {0, 2, 4}  
    \filldraw [ color=white!70!gray] (\x, \y) rectangle (\x+1, \y+1);
\foreach \x in {-3, -1, 1, 3}
  \foreach \y in {1, 3,  5}  
    \filldraw [ color=white!70!gray] (\x, \y) rectangle (\x+1, \y+1);
\foreach \y in {0, 2, 4}
  \filldraw [color=white!70!gray] (-3.5, \y) rectangle (-3, \y+1);
\foreach \y in {1, 3, 5}
  \filldraw [color=white!70!gray] (5, \y) rectangle (5.5, \y+1);
\foreach \x in {-2, 0, 2, 4}
  \filldraw [color=white!70!gray] (\x, 6) rectangle (\x+1, 6.5);
\foreach \x in {-3, -1, 1, 3}
  \filldraw [color=white!70!gray] (\x, -0.5) rectangle (\x+1, 0);

\filldraw  [color=white!70!gray] (-3.5, 6) rectangle (-3, 6.5);
\filldraw  [color=white!70!gray] (5, -0.5) rectangle (5.5, 0);

\node  at (1.5, 3.45) {\BlackKingOnWhite};
\node  at (2.5, 3.45) {\BlackKingOnWhite};
\node  at (2.5, 2.45) {\BlackKingOnWhite};
\node  at (-0.5, 3.45) {\WhiteKingOnWhite};
\node  at (0.5, 3.45) {\WhiteBishopOnWhite};

\draw [thick, color=red] (4, 6.5)--(4,6)--(3,6)--(3,5)--(2,5)--(2,4)--(1,4)--(1,3)--(2,3)--(2,2)
--(1,2)--(1,1)--(0,1)--(0,0)--(-1,0)--(-1,-0.5);
 
 \draw[ thick, ->] (-2.5, 3.5)--(-0.9,2.7);  
 \node  at (-0.5, 2.45) {\WhiteKnightOnWhite};

\node [anchor=base] at (-2.5, -1.2) {a};
\node [anchor=base] at (-1.5, -1.2) {b};
\node [anchor=base] at (-0.5, -1.2) {c};
\node [anchor=base] at (0.5, -1.2) {d};
\node [anchor=base] at (1.5, -1.2) {e};
\node [anchor=base] at (2.5, -1.2) {f};
\node [anchor=base] at (3.5, -1.2) {g};
\node [anchor=base] at (4.5, -1.2) {h};

\node at (-4, 0.5) {$2$};
\node at (-4, 1.5) {$3$};
\node at (-4, 2.5) {$4$};
\node at (-4, 3.5) {$5$};
\node at (-4, 4.5) {$6$};
\node at (-4, 5.5) {$7$};

\node at (6, 0.5) {$2$};
\node at (6, 1.5) {$3$};
\node at (6, 2.5) {$4$};
\node at (6, 3.5) {$5$};
\node at (6, 4.5) {$6$};
\node at (6, 5.5) {$7$};

\end{tikzpicture}

\caption{The knight takes action if the black king waits around the white king and bishop.}
\label{F:action} 

\end{center}

\end{figure}
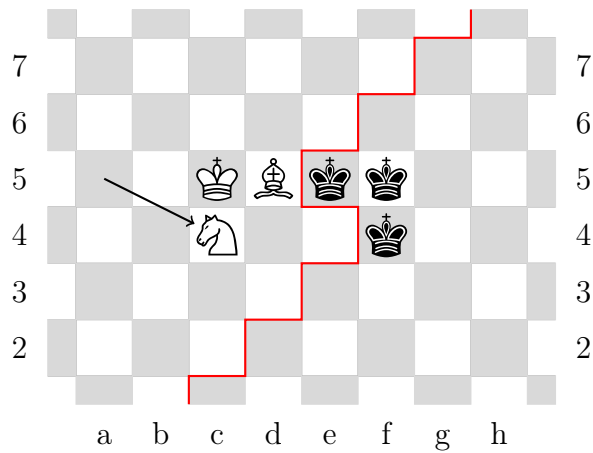

In the situation of Figure~\ref{F:action}, assuming that Black waits on the 4:th and 5:th ranks, for instance in the three squares indicated, White will maneuver the knight to c4. 
Notice that the knight on c4 doesn't get in the way of the king or bishop, and that the enclosure can be maintained whatever Black does. Since the black king now can't go to e3, the white king doesn't need access to c4. 

Suppose first that White gets to play 1.~Nc4 while Black is on the 4:th rank or higher. Black can now no longer move to the 3:rd rank without allowing White to shift the bishop to e4 and establish a smaller simple enclosure. 

If Black stays on ranks 4--6, White plays 2.~Kd6 and tries to decrease the enclosure by getting the bishop either to e4 or to e6. The only way Black can try to prevent this is by playing 2.~- Kf4/f5. But then White can force Black away and get the bishop to e4 anyway, for instance 2.~- Kf5 3.~Ke7 Kg4 4.~Kf6 Kf4 5.~Ke6, forcing Black to the g-file, and 6.~Be4 establishing the simple enclosure. 

There is one slight problem that needs to be treated separately:  
If just before the knight is about to go to c4, Black plays 1.~- Ke3, then White is forced to play 2.~Kc4. If Black now returns with 2.~- Kf4, White can't play 3.~Nc4 because the king is in the way (if Black instead continues south, White shifts the bishop to e4, making progress without the knight). But White then triangulates with 3.~Kd4 Kf5/g5 (otherwise the bishop gets to e4) 4.~Kc5, and unless Black goes to the 6:th rank, 5.~Nc4.

In some situations, other knight maneuvers can be faster and just as effective (the knight just needs to control some key dark squares), but going to c4 (or the corresponding square relative to the bishop) always works. 

\subsection*{Conclusion: White can force checkmate}
Playing White from any position, and provided we can avoid immediate stalemate or losing a piece, we can centralize our pieces and establish an enclosure in some direction. Basically we just put our bishop next to our king and notice where the black king is relative to the bishop's diagonals. Possibly we also have to let the knight run away from the black king. Then, using the maneuvers already described, we can make sure to eventually decrease the enclosure. There is only one way that this process can end: We finally constrain the king to the immediate vicinity of a mating corner, and the standard technique applies.

\section{Number of moves required for checkmate}

In endgames like rook or two bishops against a bare king, checkmate can be forced in time (number of moves) proportional to the side of the board, let's call it $n$. The Telesin checkmating procedure for KBN vs K on the other hand requires time proportional to the area $n^2$ of the board. When White reduces to a simple enclosure, Black can choose in which direction to run, and can therefore run back and forth along the barrier, eventually visiting a positive proportion of the squares within the original enclosure. 

Telesin states that a 1000 by 1000 board allows checkmate in at most 1,028,000 moves, indicating that his calculations yield a main term of $n^2$ plus some overhead. In \cite{oeis} it is conjectured that from the worst possible starting configuration that allows a forced win at all, the checkmate requires about $n^2/6+3n$ moves for even $n$, and $n^2/6 + 2n$ for odd $n$. The conjecture seems to be due to Kot\v{e}\v{s}ovec, but it's not entirely clear from \cite{Kotesovec, oeis} how strong the evidence is.     

\subsection*{Forcing along an edge anyway?}
We believe that Telesin's method is generally necessary even if we have forced Black to a safe corner, since it seems we cannot get the black king out of such a corner without also letting it escape almost to the opposite safe corner. If this is the case, then a position with the black king trapped in a safe corner is, for large $n$, essentially the worst case, since the king is at maximum distance from a mating corner and requires the maximum number of rounds of zone restriction. 

But this is still not entirely clear. As if just to make things even more confusing, it \emph{is} possible to force the black king along an edge of the board, if only we can get it away from the safe corner without misplacing our pieces! This is another (minor) point that makes the Telesin article a bit hard to understand: In the main line the black king runs immediately to an edge of the board, from where it would have been possible (in \cite[§7]{Telesin}) to shove it directly to a mating corner, checkmating faster than in the line given.   

Suppose we arrive at the position below, with Black to move and assuming that White tries to drive the black king to a distant light square corner to the right. 

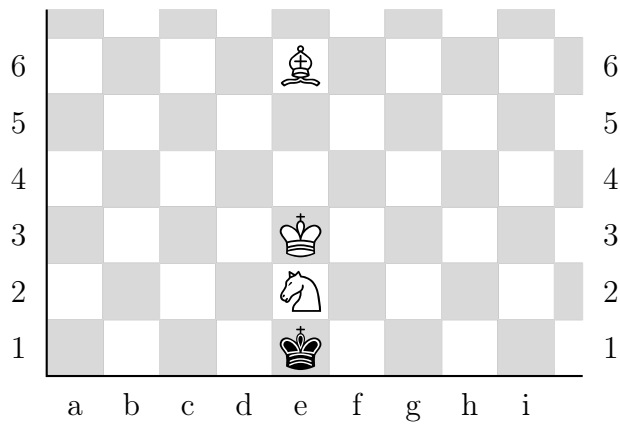
\begin{center}
\begin{tikzpicture} [scale=0.75]

\draw[step=1cm,gray,ultra thin] (-4, 0) grid (5.5,6.5); 

\foreach \x in {-4, -2, 0, 2, 4}
  \foreach \y in {0, 2, 4}  
    \filldraw [ color=white!70!gray] (\x, \y) rectangle (\x+1, \y+1);
\foreach \x in {-3, -1, 1, 3}
  \foreach \y in {1, 3,  5}  
    \filldraw [ color=white!70!gray] (\x, \y) rectangle (\x+1, \y+1);

\foreach \y in {1, 3, 5}
  \filldraw [color=white!70!gray] (5, \y) rectangle (5.5, \y+1);
\foreach \x in {-4, -2, 0, 2, 4}
  \filldraw [color=white!70!gray] (\x, 6) rectangle (\x+1, 6.5);

\draw [thick, color=black] (-4,6.5)--(-4,0)--(5.5, 0);

\node  at (0.5, 2.45) {\WhiteKingOnWhite};
\node  at (0.5, 0.45) {\BlackKingOnWhite};
\node  at (0.5, 1.45) {\WhiteKnightOnWhite};
\node  at (0.5, 5.45) {\WhiteBishopOnWhite};

\node [anchor=base] at (-3.5, -0.7) {a};
\node [anchor=base] at (-2.5, -0.7) {b};
\node [anchor=base] at (-1.5, -0.7) {c};
\node [anchor=base] at (-0.5, -0.7) {d};
\node [anchor=base] at (0.5, -0.7) {e};
\node [anchor=base] at (1.5, -0.7) {f};
\node [anchor=base] at (2.5, -0.7) {g};
\node [anchor=base] at (3.5, -0.7) {h};
\node [anchor=base] at (4.5, -0.7) {i};

\node at (-4.5, 0.5) {$1$};
\node at (-4.5, 1.5) {$2$};
\node at (-4.5, 2.5) {$3$};
\node at (-4.5, 3.5) {$4$};
\node at (-4.5, 4.5) {$5$};
\node at (-4.5, 5.5) {$6$};

\node at (6, 0.5) {$1$};
\node at (6, 1.5) {$2$};
\node at (6, 2.5) {$3$};
\node at (6, 3.5) {$4$};
\node at (6, 4.5) {$5$};
\node at (6, 5.5) {$6$};
 
\end{tikzpicture}
\end{center}

White can then force the same relative position of the pieces to occur again, but with the entire position shifted two steps to the right, so that all pieces are on the g-file. Moreover, this will happen in at most 12 moves.

If Black starts 1.~- Kf1, White responds 2.~Nf4. If on the other hand Black chooses 1.~- Kd1, we can play 2.~Bb3+ (this is what is not possible two files further left, when digging Black out of the corner!) Ke1 3.~Nf4 Kf1 4.~Be6, so that in any case we reach the following position:
  
 
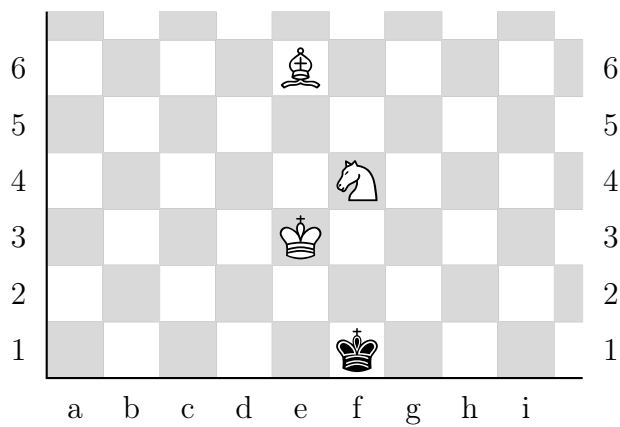
\begin{center}
\begin{tikzpicture} [scale=0.75]

\draw[step=1cm,gray,ultra thin] (-4, 0) grid (5.5,6.5); 

\foreach \x in {-4, -2, 0, 2, 4}
  \foreach \y in {0, 2, 4}  
    \filldraw [ color=white!70!gray] (\x, \y) rectangle (\x+1, \y+1);
\foreach \x in {-3, -1, 1, 3}
  \foreach \y in {1, 3,  5}  
    \filldraw [ color=white!70!gray] (\x, \y) rectangle (\x+1, \y+1);

\foreach \y in {1, 3, 5}
  \filldraw [color=white!70!gray] (5, \y) rectangle (5.5, \y+1);
\foreach \x in {-4, -2, 0, 2, 4}
  \filldraw [color=white!70!gray] (\x, 6) rectangle (\x+1, 6.5);

\draw [thick, color=black] (-4,6.5)--(-4,0)--(5.5, 0);

\node  at (0.5, 2.45) {\WhiteKingOnWhite};
\node  at (1.5, 0.45) {\BlackKingOnWhite};
\node  at (1.5, 3.45) {\WhiteKnightOnWhite};
\node  at (0.5, 5.45) {\WhiteBishopOnWhite};

\node [anchor=base] at (-3.5, -0.7) {a};
\node [anchor=base] at (-2.5, -0.7) {b};
\node [anchor=base] at (-1.5, -0.7) {c};
\node [anchor=base] at (-0.5, -0.7) {d};
\node [anchor=base] at (0.5, -0.7) {e};
\node [anchor=base] at (1.5, -0.7) {f};
\node [anchor=base] at (2.5, -0.7) {g};
\node [anchor=base] at (3.5, -0.7) {h};
\node [anchor=base] at (4.5, -0.7) {i};

\node at (-4.5, 0.5) {$1$};
\node at (-4.5, 1.5) {$2$};
\node at (-4.5, 2.5) {$3$};
\node at (-4.5, 3.5) {$4$};
\node at (-4.5, 4.5) {$5$};
\node at (-4.5, 5.5) {$6$};

\node at (6, 0.5) {$1$};
\node at (6, 1.5) {$2$};
\node at (6, 2.5) {$3$};
\node at (6, 3.5) {$4$};
\node at (6, 4.5) {$5$};
\node at (6, 5.5) {$6$};
 
\end{tikzpicture}
\end{center}
Now on 4.~- Kg1 White plays 5.~Kf3 leading to the position below, while if Black tries to go back, we have 4.~Ke1 5.~Nd3+! Here on 5.~- Kd1 Black gets checkmated immediately with 6.~Bb3\# which would not be possible two files further left. Therefore 5.~- Kf1 is forced, and after 6.~Kf3 Kg1 7.~Nf4 we reach the same position:

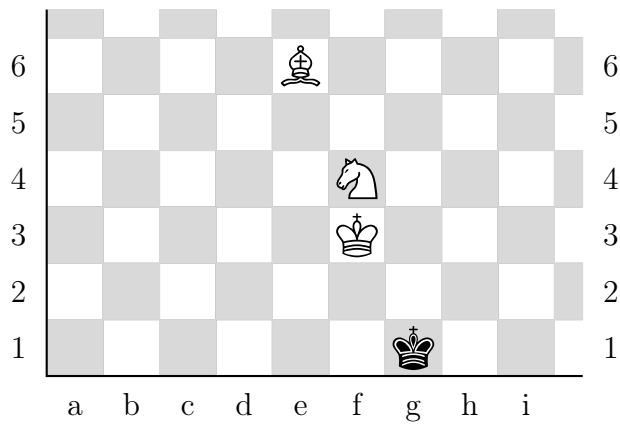
\begin{center}
\begin{tikzpicture} [scale=0.75]

\draw[step=1cm,gray,ultra thin] (-4, 0) grid (5.5,6.5); 

\foreach \x in {-4, -2, 0, 2, 4}
  \foreach \y in {0, 2, 4}  
    \filldraw [ color=white!70!gray] (\x, \y) rectangle (\x+1, \y+1);
\foreach \x in {-3, -1, 1, 3}
  \foreach \y in {1, 3,  5}  
    \filldraw [ color=white!70!gray] (\x, \y) rectangle (\x+1, \y+1);
\foreach \y in {1, 3, 5}
  \filldraw [color=white!70!gray] (5, \y) rectangle (5.5, \y+1);
\foreach \x in {-4, -2, 0, 2, 4}
  \filldraw [color=white!70!gray] (\x, 6) rectangle (\x+1, 6.5);

\draw [thick, color=black] (-4,6.5)--(-4,0)--(5.5, 0);

\node  at (1.5, 2.45) {\WhiteKingOnWhite};
\node  at (2.5, 0.45) {\BlackKingOnWhite};
\node  at (1.5, 3.45) {\WhiteKnightOnWhite};
\node  at (0.5, 5.45) {\WhiteBishopOnWhite};

\node [anchor=base] at (-3.5, -0.7) {a};
\node [anchor=base] at (-2.5, -0.7) {b};
\node [anchor=base] at (-1.5, -0.7) {c};
\node [anchor=base] at (-0.5, -0.7) {d};
\node [anchor=base] at (0.5, -0.7) {e};
\node [anchor=base] at (1.5, -0.7) {f};
\node [anchor=base] at (2.5, -0.7) {g};
\node [anchor=base] at (3.5, -0.7) {h};
\node [anchor=base] at (4.5, -0.7) {i};

\node at (-4.5, 0.5) {$1$};
\node at (-4.5, 1.5) {$2$};
\node at (-4.5, 2.5) {$3$};
\node at (-4.5, 3.5) {$4$};
\node at (-4.5, 4.5) {$5$};
\node at (-4.5, 5.5) {$6$};

\node at (6, 0.5) {$1$};
\node at (6, 1.5) {$2$};
\node at (6, 2.5) {$3$};
\node at (6, 3.5) {$4$};
\node at (6, 4.5) {$5$};
\node at (6, 5.5) {$6$};
 
\end{tikzpicture}
\end{center}

Next White plays 8.~Ng2 no matter where Black goes. Since the bishop already controls h3 and i2, Black cannot get further than to h2. White can then maneuver the king to g3 and the bishop to g6 while making sure Black never escapes to e2 or i2. No later than after White's move 13 we arrive at:

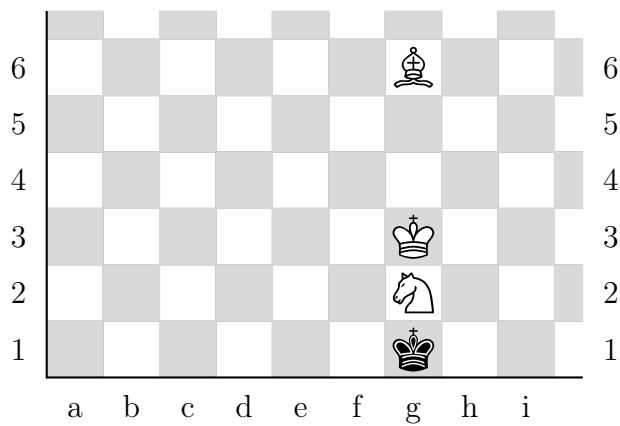
\begin{center}
\begin{tikzpicture} [scale=0.75]

\draw[step=1cm, gray, ultra thin] (-4, 0) grid (5.5,6.5); 

\foreach \x in {-4, -2, 0, 2, 4}
  \foreach \y in {0, 2, 4}  
    \filldraw [ color=white!70!gray] (\x, \y) rectangle (\x+1, \y+1);
\foreach \x in {-3, -1, 1, 3}
  \foreach \y in {1, 3,  5}  
    \filldraw [ color=white!70!gray] (\x, \y) rectangle (\x+1, \y+1);

\foreach \y in {1, 3, 5}
  \filldraw [color=white!70!gray] (5, \y) rectangle (5.5, \y+1);
\foreach \x in {-4, -2, 0, 2, 4}
  \filldraw [color=white!70!gray] (\x, 6) rectangle (\x+1, 6.5);

\draw [thick, color=black] (-4,6.5)--(-4,0)--(5.5, 0);

\node  at (2.5, 2.45) {\WhiteKingOnWhite};
\node  at (2.5, 0.45) {\BlackKingOnWhite};
\node  at (2.5, 1.45) {\WhiteKnightOnWhite};
\node  at (2.5, 5.45) {\WhiteBishopOnWhite};

\node [anchor=base] at (-3.5, -0.7) {a};
\node [anchor=base] at (-2.5, -0.7) {b};
\node [anchor=base] at (-1.5, -0.7) {c};
\node [anchor=base] at (-0.5, -0.7) {d};
\node [anchor=base] at (0.5, -0.7) {e};
\node [anchor=base] at (1.5, -0.7) {f};
\node [anchor=base] at (2.5, -0.7) {g};
\node [anchor=base] at (3.5, -0.7) {h};
\node [anchor=base] at (4.5, -0.7) {i};

\node at (-4.5, 0.5) {$1$};
\node at (-4.5, 1.5) {$2$};
\node at (-4.5, 2.5) {$3$};
\node at (-4.5, 3.5) {$4$};
\node at (-4.5, 4.5) {$5$};
\node at (-4.5, 5.5) {$6$};

\node at (6, 0.5) {$1$};
\node at (6, 1.5) {$2$};
\node at (6, 2.5) {$3$};
\node at (6, 3.5) {$4$};
\node at (6, 4.5) {$5$};
\node at (6, 5.5) {$6$};
 
\end{tikzpicture}
\end{center}

The pattern then repeats, and when we approach the south-east corner, the standard checkmating method applies. Therefore as soon as we reach a position with the black king trapped in a controlled way along an edge and not too close to a safe corner, we can obtain checkmate in linear time. 

We have to leave the question of the ``time complexity'' of the KBN vs K checkmate open.
An interesting challenge would be to describe a strategy for Black that survives from some ``stable'' position for $c\cdot n^2$ moves, for some fixed $c>0$. But perhaps White can restrict the black king in some slightly more efficient way, nudging it towards a mating corner at a small but constant speed, so that checkmate is forced in linear time after all.

\end{document}